\documentclass[final,1p,times]{elsarticle}

\usepackage{amssymb}
 \usepackage{amsthm}
\usepackage{amscd}
\usepackage{amsmath}
\usepackage{amsfonts}
\usepackage{amssymb}
\usepackage{graphicx}
\newtheorem{theorem}{Theorem}

\newtheorem{lemma}[theorem]{Lemma}

\usepackage{mathrsfs}
\usepackage{titletoc}

\newcommand{\la}{\Delta}
\def\na{\nabla}

\renewcommand{\d}{\delta}

\newcommand{\ra}{\rightarrow}

\newcommand{\f}{\frac}

\def\e{\epsilon}

\def\wan{\widetilde}

\newcommand{\be}{\begin{equation}}
\renewcommand{\ra}{\rightarrow}
\newcommand{\ee}{\end{equation}}
\newcommand{\bea}{\begin{eqnarray}}
\newcommand{\eea}{\end{eqnarray}}
\newcommand{\bna}{\begin{eqnarray*}}
\newcommand{\ena}{\end{eqnarray*}}

\renewcommand{\O}{\Omega}
\renewcommand{\le}{\left}
\newcommand{\ri}{\right}

\journal{***}

\begin{document}

\begin{frontmatter}

\title{A remark on a result of Ding-Jost-Li-Wang}

\author{Yunyan Yang}
 \ead{ yunyanyang@ruc.edu.cn}
 \author{Xiaobao Zhu}
  \ead{zhuxiaobao@ruc.edu.cn}
\address{ Department of Mathematics,
Renmin University of China, Beijing 100872, P. R. China}

\begin{abstract}

Let $(M,g)$ be a compact Riemannian surface without boundary, $W^{1,2}(M)$ be the usual Sobolev space,
$J: W^{1,2}(M)\ra \mathbb{R}$ be the functional defined by
$$J(u)=\f{1}{2}\int_M|\nabla u|^2dv_g+8\pi \int_M udv_g-8\pi\log\int_Mhe^udv_g,$$
where $h$ is a positive smooth function on $M$. In an inspiring work (Asian J. Math., vol. 1, pp. 230-248, 1997),
Ding, Jost, Li and Wang obtained a sufficient condition under which $J$ achieves its minimum.
In this note, we prove that if the smooth function $h$ satisfies $h\geq 0$ and $h\not\equiv 0$, then
the above result still holds. Our method is to exclude blow-up points on the zero set of $h$.
\end{abstract}

\begin{keyword}
Kazdan-Warner problem, Trudinger-Moser inequality, blow-up analysis
\MSC[2010] 46E35
\end{keyword}

\end{frontmatter}

\titlecontents{section}[0mm]
                       {\vspace{.2\baselineskip}}
                       {\thecontentslabel~\hspace{.5em}}
                        {}
                        {\dotfill\contentspage[{\makebox[0pt][r]{\thecontentspage}}]}
\titlecontents{subsection}[3mm]
                       {\vspace{.2\baselineskip}}
                       {\thecontentslabel~\hspace{.5em}}
                        {}
                       {\dotfill\contentspage[{\makebox[0pt][r]{\thecontentspage}}]}

\setcounter{tocdepth}{2}


\section{Introduction}

\subsection{A review of a result of Ding-Jost-Li-Wang}

Let $(M,g)$ be a compact Riemannian surface without boundary, $h(x)$ be a smooth function on $M$.
Suppose the area of $M$ is equal to $1$. In \cite{DJLW97}, Ding, Jost, Li and Wang
studied one of the Kazdan-Warner problems \cite{KW74}, existence of solutions to the equation
\be\label{1}\Delta u=8\pi-8\pi h e^{u}\quad{\rm on}\quad M,\ee
where $\Delta$ denotes the usual Laplacian on $M$. Their method of solving this problem is calculus of variations.
Namely, they tried to minimize the functional
$$J(u)=\f{1}{2}\int_M|\nabla u|^2dv_g+8\pi\int_M udv_g-8\pi \log\int_M he^udv_g$$
in the function space
$$H_1=\le\{u\in W^{1,2}(M)\le|\int_Mhe^udv_g=1\ri.\ri\}.$$
To show that $J$ is bounded from below, they considered subcritical functional
$$J_\epsilon(u)=\f{1}{2}\int_M|\nabla u|^2dv_g+(8\pi-\epsilon)\int_M u dv_g-(8\pi-\epsilon) \log\int_M he^udv_g$$
for any $\epsilon>0$. It is not difficult to check that $J$ achieves its minimum in the function space $H_1$. The minimizer
$u_\epsilon$ satisfies the Euler-Lagrange equation
\be\label{e-l}\Delta u_\epsilon=(8\pi-\epsilon)-(8\pi-\epsilon)he^{u_\epsilon}.\ee
By the elliptic regularity theory, $u_\epsilon\in C^2(M)\cap H_1$. Set
$\lambda_\epsilon=\max_{M} u_\epsilon$. There are two possibilities:
$(i)$ there exists some constant $C$ such that $\lambda_\epsilon\leq C$ for all $\epsilon>0$;
$(ii)$ $\lambda_\epsilon\ra +\infty$ as $\epsilon\ra 0$, i.e., $u_\epsilon$ blows up.
In case $(i)$, by applying  elliptic estimates to (\ref{e-l}), they proved that $u_\epsilon\ra u_0$ in $C^1(M)$
and that $u_0$ is a minimizer of $J$ in $H_1$. In case $(ii)$, they assumed that $u_\epsilon(x_\epsilon)=\lambda_\epsilon$
and $x_\epsilon\ra p\in M$. Here and throughout this note, sequence and subsequence are not distinguished.
Choose a local normal coordinate system around $p$. Let $\lambda_\epsilon^\ast=e^{\lambda_\epsilon/2}$ and
$\varphi_\epsilon(x)=u_\epsilon(x_\epsilon+x/{\lambda_\epsilon^\ast})-\lambda_\epsilon$.
By elliptic estimates, they showed that $\varphi_\epsilon\ra \varphi$ in $C^1_{\rm loc}(\mathbb{R}^2)$, where $\varphi$ satisfies
\be\label{bubble}
\le\{\begin{array}{lll}
\Delta_{\mathbb{R}^2}\varphi(x)=-8\pi h(p)e^{\varphi(x)}\\
[1.5ex]\varphi(0)=0=\sup_{\mathbb{R}^2}\varphi\\[1.5ex]
\int_{\mathbb{R}^2}h(p)e^{\varphi(x)}dx\leq 1,
\end{array}\ri.\ee
where $\Delta_{\mathbb{R}^2}$ denotes the standard Laplacian on $\mathbb{R}^2$.
Suppose that
$h(x)>0$  for all $x\in M$.
In particular, $h(p)>0$. It then follows from a classification theorem of Chen-Li \cite{CL91} that
$$\varphi(x)=-2\log(1+\pi h(p)|x|^2).$$
Moreover, they proved that $u_\epsilon(x)-\overline{u}_\epsilon(x)\ra G(x,p)$ weakly in
$W^{1,q}(M)$ for any $1<q<2$ and in $C^2_{\rm loc}(M\setminus\{p\})$,
where $\overline{u}_\epsilon=\int_M u_\epsilon dv_g$ and $G$ is a Green function satisfying
\begin{align*}
\begin{cases}
      & \la G=8\pi-8\pi\d_p, \\
      & \int_M Gdv_g=0.
\end{cases}
\end{align*}
In a normal coordinate system around $p$ one can assume that
\begin{align}
\label{green-1}
G(x,p)=-4\log r+A(p)+b_1x_1+b_2x_2+c_1x_1^2+2c_2x_1x_2+c_3x_2^2+O(r^3),
\end{align}
where $r(x)=dist(x,p)$.

Furthermore, using the maximum principle, they analyzed the neck energy of $u_\epsilon$.
Combining all the above analysis, they concluded that if $u_\epsilon$ blows up, then
\be\label{low-bd}\inf_{u\in H_1}J(u)\geq C_0=-8\pi-8\pi\log\pi-4\pi\max_{p\in M}(A(p)+2\log h(p)).\ee
In other words, if (\ref{low-bd}) does not hold, then $u_\epsilon$ would not blow up and thus $C^1$-converges to
a minimizer of $J$.
Then the proof of existence result for (\ref{1}) was reduced to constructing a sequence of functions $\phi_\epsilon$ satisfying
$J(\phi_\epsilon)< C_0$
for sufficiently small $\epsilon>0$. They proved that this is true provided that
\be\label{cond}\la h(p_0)+2(b_1(p_0)k_1(p_0)+b_2(p_0)k_2(p_0))>-(8\pi+(b_1^2(p_0)+b_2^2(p_0))-2K(p_0))h(p_0),\ee
where $K(x)$ is the Gaussian curvature on $M$, $\nabla h(p_0)=(k_1(p_0),k_2(p_0))$ in the normal coordinate system,
$p_0$ is the maximum point of the function $A(q)+2\log h(q)$, and $b_1(p_0)$, $b_2(p_0)$ are constants in the expression of
the Green function $G$, namely (\ref{green-1}).

Summarizing, they obtained the following result (\cite{DJLW97}, Theorem 1.2):
If $h$ is a smooth and strictly positive function on $M$ and
the hypothesis (\ref{cond}) is satisfied, then the equation (\ref{1}) has a smooth solution.

\subsection{An improvement}

Checking the proof of (\cite{DJLW97}, Theorem 1.2), we find that the hypothesis $h(x)>0$ for all $x\in M$ is only needed in
solving the bubble $\varphi$ from (\ref{bubble}), where $p$ is the blow-up point. Now suppose that $h\geq 0$
and $h\not\equiv 0$ on $M$. If $u_\epsilon$ blows up at $p\in M$ with $h(p)>0$, then we conclude that (\cite{DJLW97}, Theorem 1.2)
still holds. Precisely, we have the following:

\begin{theorem}\label{mainthm}
Let $(M,g)$ be a compact Riemann surface, $K(x)$ be its Gauss curvature, and
$h: M\ra \mathbb{R}$ be a smooth function such that $h(x)\geq 0$ and $h\not\equiv 0$.
Suppose that $A(x)+2\log h(x)$ attains its maximum at $p^\ast\in M$. Let $b_1(p^\ast)$ and $b_2(p^\ast)$ be two constants
defined as in (\ref{green-1}), and write $\na h(p^\ast)=(k_1(p^\ast),k_2(p^\ast))$ in the normal coordinate system. If
(\ref{cond}) holds at $p^\ast$,
then (\ref{1}) has a smooth solution.
\end{theorem}

The method we exclude the possibility that $u_\epsilon$ blows up at some zero point of $h$ is based on a concentration lemma
of Ding-Jost-Li-Wang \cite{DJLW98}.

Before ending the introduction, we mention two related works also about the existence results of mean field equations 
at the critical parameter. Lin-Wang \cite{LW2010} obtained a sufficient and necessary condition under which the mean 
field equation on a flat  torus at the critical parameter has a solution. Bartolucci-Lin \cite{BL2014} considered the existence 
of mean field equations at the critical parameter on multiply connected domains in $\mathbb{R}^2$.

\section{Preliminaries}
\subsection{The Trudinger-Moser inequality}

We state a Trudinger-Moser inequality, which plays an important role in the Kazdan-Warner problem (\ref{1}).
\begin{lemma}\label{MT}
(\cite{Fontana,DJLW97}) Let $(M,g)$ be a compact Riemannian surface. For any $u\in W^{1,2}(M)$ with $\int_M udv_g=0$
one has
 \begin{align*}
  \int_M e^u dv_g\leq C_M e^{\frac{1}{16\pi}\|\na u\|_2^2},
 \end{align*}
where $C_M$ is a positive constant depending only on $(M,g)$.
\end{lemma}

\subsection{The concentration lemma}

The key tool we use to solve Theorem \ref{mainthm} is the following concentration lemma due to
Ding-Jost-Li-Wang.
\begin{lemma}\label{prop1}(\cite{DJLW98})
Let $(M,g)$ be a compact Riemannian surface with volume $1$. Given a sequence of functions $v_\epsilon\in W^{1,2}(M)$ with
$\int_M e^{v_\epsilon}dv_g=1$ and $\int_M|\na v_\epsilon|^{2}dv_g+16\pi\int_M v_\epsilon dv_g\leq C$, then either\\
$(i)$ there is a constant $C_1>0$ such that $\int_M|\na v_\epsilon|^{2}dv_g\leq C_1$; or\\
$(ii)$ $v_\epsilon$ concentrates at a point $p\in M$, i.e., for any $r>0$,
$$\lim_{\epsilon\ra 0}\int_{B_r(p)}e^{v_\epsilon}dv_g=1.$$
\end{lemma}

\section{Proof of Theorem \ref{mainthm}}

{\it Proof of Theorem \ref{mainthm}}.
Let $\e>0$, we consider the functional
\bna
J_\e(u)=\frac{1}{2}\int_M |\na u|^2dv_g+(8\pi-\e)\int_M u dv_g-(8\pi-\e)\log\int_M he^udv_g
\ena
in the function space $H_1=\le\{u\in W^{1,2}(M)\le|\int_M he^udv_g=1\ri.\ri\}$. By Lemma \ref{MT},
$J_\epsilon$ is coercive and bounded from below, so it attains its infimum
in $H_1$ at some $u_\epsilon\in H_1$. Clearly,
$u_\epsilon$ satisfies the Euler-Lagrange equation (\ref{e-l}), i.e.
\bna
\la u_\e=(8\pi-\e)-(8\pi-\e)he^{u_\e}.
\ena
Denote $\wan{H}=\le\{u\in W^{1,2}(M)\le|\int_M udv_g=0\ri.\ri\}$. For any $u\in W^{1,2}(M)$, we have
$u-\overline{u}\in \wan{H}$ and $u-\log\int_Mhe^udv_g \in H_1$, where $\overline{u}=\int_Mudv_g$.
Thus we have
$$\inf_{u\in H_1}J_\e(u)=\inf_{u\in\wan{H}}J_\e(u)=\inf_{u\in W^{1,2}(M)}J_\e(u),~~\inf_{u\in H_1}J(u)=\inf_{u\in\wan{H}}J(u)=
\inf_{u\in W^{1,2}(M)}J(u),$$
since $J_\e(u+c)=J_\e(u)$ and $J(u+c)=J(u)$ for any constant $c$.
 Let $w_\e=u_\e-\overline{u}_\epsilon$, then  $w_\epsilon\in \widetilde{H}$ attains the infimum of
 $J_\e$ in $\widetilde{H}$.
There are two possibilities:\\

\textbf{Case 1.} $||\na u_\e||_2\leq C$.\\

In this case, $\|\na w_\e\|_2\leq C$. Without loss of generality, we assume
\bna
&&w_\e\rightharpoonup w_0\quad{\rm weakly\,\,in}\quad W^{1,2}(M)\\
&&w_\e\ra w_0\quad {\rm strongly\,\,in}\quad L^p(M),\,\,\forall p\geq1.
\ena
This together with the Trudinger-Moser inequality (Lemma \ref{MT}) and the H\"{o}lder inequality
leads to
\begin{align*}
\label{}
\int_Mh\le(e^{w_\e}-e^{w_0}\ri)dv_g
&=\int_M h \int_0^1 \frac{d}{dt}e^{w_0+t(w_\e-w_0)}dtdv_g\nonumber\\
&=\int_0^1 \int_M h e^{w_0+t(w_\e-w_0)}(w_\e-w_0)dv_g dt\nonumber\\
&\ra0~~as~~\e\ra0.
\end{align*}
It follows that
\begin{align*}
\label{}
\inf_{u\in \wan{H}}J(u)\geq\liminf_{\e\ra0}\inf_{u\in\wan{H}}J_\e(u)=\liminf_{\e\ra0}J_\e(w_\e)\geq J(w_0).
\end{align*}
That is to say, $w_0\in \wan{H}$ attains the infimum of $J$ in $\wan{H}$. Let $u_0=w_0-\log\int_M he^{w_0}dv_g$. Hence $u_0\in H_1$
attains the infimum of $J$ in $H_1$ and satisfies
$$\la u_0=8\pi-8\pi he^{u_0}.$$
The elliptic regularity theory implies that $u_0\in C^2(M)$. The proof of Theorem \ref{mainthm} terminates in this case.\\

\textbf{Case 2.} $||\na u_\e||_2\ra+\infty$. \\

By the elliptic regularity theory, in view of (\ref{e-l}), we have $u_\epsilon\in C^2(M)\cap H_1$. Let
$\lambda_\epsilon=\max_{M} u_\e=u_\e(x_\e)$. In this case, there holds $\lambda_\epsilon\ra+\infty$.
For otherwise, we have that $\|\nabla u_\epsilon\|_2$ is bounded by applying elliptic estimates to (\ref{e-l}).
It should be remarked that $\|\nabla u_\epsilon\|_2\ra +\infty$ if and only if $\lambda_\epsilon\ra+\infty$.
Let $\O\subset M$ be a domain. If
$\int_\O he^{u_\e}dv_g\leqslant\frac{1}{2}-\d$
for some $\d\in(0,\frac{1}{2})$, then (\cite{DJLW97}, Lemma 2.8) implies that
\be\label{Fact3}||u_\e-\overline{u}_\e||_{L^{\infty}(\O_0)}\leq C(\O_0,\O),\quad\forall \O_0\subset\subset\O.\ee

Suppose $x_\e\ra p$ as $\e\ra0$. As we explained in the introduction, to prove Theorem \ref{mainthm}, it
 suffices to prove that $h(p)>0$. For this purpose, we set
$v_\e=u_\e-\log\int_M e^{u_\e}dv_g$. Then $\int_M e^{v_\e}dv_g=1$. Note that
$J_\e(v_\e)=J_\e(u_\e)$ and that
$$J_\e(u_\e)=\inf_{H_1}J_\e(u)\leq J_\e\le(-\log\int_M h dv_g\ri)=-(8\pi-\e)\log\int_M h dv_g.$$
By Jensen's inequality, we have
\be\label{mean}\int_M v_\e dv_g\leq\log\int_M e^{v_\e}dv_g=0.\ee
Combining the above two estimates, we obtain
\begin{align}
\int_M|\na v_\e|^2dv_g+16\pi\int_M v_\e dv_g\leq& 2\left(J_\e(v_\e)+\e\int_Mv_\e dv_g+(8\pi-\e)\log\int_Mhe^{v_\e}dv_g\right)\nonumber\\
\leq&2\left(J_\e(u_\e)+\e\int_Mv_\e dv_g+(8\pi-\e)\log\max_M h \right)\nonumber\\
\leq&(16\pi-2\e)\log\frac{\max_M{h}}{\int_M h dv_g}\nonumber\\
\leq&C.\nonumber
\end{align}
Clearly, $(ii)$ of Lemma \ref{prop1} holds in this case. Hence there exists some $p^\prime\in M$ such that
$v_\e$ concentrates at $p'$, namely,
\be\label{conc}\lim_{\e\ra0}\int_{B_r(p')}e^{v_\e}dv_g=1,\quad\forall r>0.\ee
We first claim that
\be\label{claim}
h(p')>0.\ee
To see this, in view of (\ref{conc}), we calculate
\begin{align}
\label{cal}
\frac{ \int_M he^{u_\e}dv_g}{\int_M e^{u_\e}dv_g}=\int_M he^{v_\e}dv_g&=\int_{B_r(p')}he^{v_\e}dv_g+
\int_{M\setminus B_r(p')}he^{v_\e}dv_g\nonumber\\
                          &=(h(p')+o_r(1))\int_{B_r(p')}e^{v_\e}dv_g+o_\e(1)\nonumber\\
                          &=(h(p')+o_r(1))(1+o_\e(1))+o_\e(1)\nonumber\\
                          &=h(p')+o_r(1)+o_\e(1).
\end{align}
Because the left hand side of (\ref{cal}) does not depend on $r$,  we have by passing to the limit
$r\ra 0$,
\begin{align}
\label{pp}
\frac{ \int_M he^{u_\e}dv_g}{\int_M e^{u_\e}dv_g}=h(p')+o_\e(1).
\end{align}
Also we have
\begin{align}
\label{17}
-(8\pi-\e)\log\int_M hdv_g= J_\e\le(-\log\int_M hdv_g\ri)\geq\inf_{H_1}J_\e(u)= J_\e(u_\e)
\end{align}
and
\begin{align}
\label{18}
J_\e(u_\e)=&\frac{1}{2}\int_M |\na u_\e|^2dv_g+(8\pi-\e)\int_M u_\e dv_g-(8\pi-\e)\log\int_M he^{u_\e}dv_g\nonumber\\
                =&\frac{1}{2}\int_M |\na u_\e|^2dv_g+(8\pi-\e)\int_M u_\e dv_g-(8\pi-\e)\log\le(h(p')+o_\e(1)\ri)\nonumber\\
                    &-(8\pi-\e)\log\int_Me^{u_\e}dv_g\nonumber\\
                 \geq&\frac{1}{2}\int_M |\na u_\e|^2dv_g+(8\pi-\e)\int_M u_\e dv_g-(8\pi-\e)\log\le(h(p')+o_\e(1)\ri)\nonumber\\
                    &-(8\pi-\e)\log\le(C e^{\frac{1}{16\pi}\int_M |\na u_\e|^2dv_g+\int_M u_\e dv_g}\ri)\nonumber\\
                 \geq&-(8\pi-\e)\log\le(C(h(p')+o_\e(1))\ri).
\end{align}
Combining (\ref{17}) and (\ref{18}), we obtain
\bna
-(8\pi-\e)\log\int_M h dv_g\geq-(8\pi-\e)\log\le(C(h(p')+o_\e(1))\ri),
\ena
and whence
$$\log\int_M h dv_g\leq\log\le(Ch(p')\ri).$$
This immediately leads to (\ref{claim}).

We next claim that
\be\label{point}p^\prime=p.\ee
In view of $(ii)$ of Lemma \ref{prop1}, there holds $\int_\O e^{v_\e}dv_g\ra 0$,
$\forall \O\subset\subset M\setminus\{p'\}$. Noting that $u_\epsilon\in H_1$, from (\ref{pp}) we have
\begin{align}
\label{20}
\int_\O he^{u_\e}dv_g=\int_\O he^{v_\e}dv_g\int_M e^{u_\e}dv_g\leq \frac{\max_{M}h}{h(p')+o_\e(1)}
\int_\O e^{v_\e}dv_g\ra 0
\end{align}
as $\epsilon\ra 0$.
Combining (\ref{Fact3}) and (\ref{20}),  we obtain $||u_\e-\overline{u}_\e||_{L^{\infty}(\O)}\leq C$
and thus $||v_\e-\overline{v}_\e||_{L^{\infty}(\O)}\leq C$
for any $\O\subset\subset M\setminus\{p'\}$. This together with (\ref{mean}) implies that $v_\e(x)\leq C$
for all $x\in \O\subset\subset M\setminus\{p'\}$.

It follows from (\ref{pp}) and $u_\epsilon\in H_1$ that
$$\int_M e^{u_\e}dv_g=\frac{1}{h(p')+o_\e(1)}<\frac{2}{h(p')}$$
for sufficiently small $\e>0$.
Suppose $p'\neq p$. Recalling that $\lambda_\epsilon=u_\epsilon(x_\epsilon)=\max_{M}u_\epsilon$ and
 $x_\epsilon\ra p$, we find a domain $\O$ such that $x_\e\in\O \subset\subset M\setminus\{p'\}$. Hence
\bna
u_\e(x_\e)\leq v_\e(x_\e)+\log\int_M e^{u_\e}dv_g\leq C,
\ena
which contradicts $\lambda_\epsilon\ra+\infty$ and concludes our claim (\ref{point}). Combining (\ref{claim})
and (\ref{point}), we obtain $h(p)>0$. Since the remaining part of the proof of Theorem \ref{mainthm} is
completely analogous to that of (\cite{DJLW97}, Theorem 1.2), we omit the details here. $\hfill\Box$\\

{\bf Acknowledgements}. The authors are supported by the National Science Foundation of China (Grant Nos. 
11171347 and 11401575). We thank the referee for his/her careful review and very helpful suggestions.

\end{document}